\documentclass [11pt,oneside]{article}

\usepackage{amssymb}
\usepackage{amsfonts}
\usepackage{amsmath}
\usepackage{amsthm}
\usepackage{epsfig}

\newtheorem{lemma}{Lemma}[section]
\newtheorem{teo}[lemma]{Theorem}
\newtheorem{rem}[lemma]{Remark}
\newtheorem{prop}[lemma]{Proposition}
\newtheorem{cor}[lemma]{Corollary}

\newtheorem{example}[lemma]{Example}

\newcommand{\matN} {\ensuremath {\mathbb{N}}}
\newcommand{\matR} {\ensuremath {\mathbb{R}}}

\newcommand{\matS} {\ensuremath {\mathbb{S}}}
\newcommand{\matD} {\ensuremath {\mathbb{D}}}

\newcommand{\calS} {\ensuremath {\mathcal{S}}}

\newcommand{\calF} {\ensuremath {\mathcal{F}}}

\newcommand{\ordtred}{\matD^{\, 3}_{\rm o}}
\newcommand{\cyctred}{\matD^{\, 3}_{\rm c}}
\newcommand{\vertred}{\matD^{\, 3}_{\rm v}}
\newcommand{\orddues}{\matS^2_{\rm o}}
\newcommand{\cycdues}{\matS^2_{\rm c}}
\newcommand{\verdues}{\matS^2_{\rm v}}
\newcommand{\ordtres}{\matS^3_{\rm o}}
\newcommand{\cyctres}{\matS^3_{\rm c}}
\newcommand{\vertres}{\matS^3_{\rm v}}

\newcommand{\nota} [1] {\caption{\footnotesize{#1}}}

\newfont{\Got}{eufm10 scaled 1200}

\font\titsc=cmcsc10 scaled 1200

\newcommand{\mettifig}[1]{\epsfig{file=#1}}

\newcommand{\timtil}{\begin{picture}(12,12)\put(2,0){$\times$}\put(2,4.5){$\sim$}\end{picture}}

\newcommand{\hatY}{\widehat{Y}}

\author{Carlo \titsc{Petronio}\thanks{Research supported
by the INTAS Project ``CalcoMet-GT'' 03-51-3663}}

\title{Spherical splitting of 3-orbifolds}

\begin{document}

\maketitle

\begin{abstract}
    \noindent
    The famous Haken-Kneser-Milnor theorem states that every $3$-manifold
    can be expressed in a unique way as a connected sum of prime $3$-manifolds.
    The analogous statement for $3$-orbifolds has been part of the folklore for
    several years, and it was commonly believed that slight variations on the
    argument used for manifolds would be sufficient to establish it. We demonstrate
    in this paper that this is not the case, proving that the
    apparently natural notion of ``essential'' system of spherical $2$-orbifolds
    is not adequate in this context. We also
    show that the statement itself
    of the theorem must be given in a substantially different way. We then prove
    the theorem in full detail, using a certain notion of ``efficient
    splitting system.''
  \vspace{4pt}

\noindent MSC (2000): 57M99 (primary), 57M50, 57M60 (secondary).
\end{abstract}

\section*{Introduction}
Since the seminal work of Thurston~\cite{thurston:notes}, the concept of orbifold
has become a central one in 3-dimensional geometric topology, because it embodies
the notion of a 3-manifold and that of a finite group action on that 3-manifold.
As a matter of fact, Thurston himself indicated how to prove that, under some
natural topological constraints, a 3-orbifold with non-empty singular locus
carries one of the 8 geometries, so it is a quotient of a geometric 3-manifold
under a finite isometric action (see~\cite{BLP} and~\cite{CHK} for modern detailed
accounts of the argument).

In the case of manifolds (closed orientable ones, say) the two basic obstructions to the
existence of a geometric structure are related to the existence of essential
spheres and of essential tori. So one restricts first to irreducible manifolds
and next to Seifert and atoroidal ones. For this restricted class of manifolds then one
has Thurston's \emph{geometrization conjecture} (or \emph{theorem}, according
to the recent work of Perelman~\cite{Perelman}). This result can be regarded as a
\emph{uniformization theorem} for 3-manifolds in view of the
Haken-Kneser-Milnor theorem (canonical splitting of a 3-manifold along spheres
into irreducible 3-manifolds) and the Jaco-Shalen-Johansson theorem
(canonical splitting of an irreducible 3-manifold along tori
into Seifert and atoroidal 3-manifolds).

The case of 3-orbifolds is similar, where the first basic obstruction to the
existence of a geometric structure is removed by restricting to (suitably defined)
\emph{irreducible} 3-orbifolds. The natural expectation would then be to have
an analogue of the Haken-Kneser-Milnor theorem, asserting that each 3-orbifold
canonically splits along spherical 2-orbifolds into irreducible 3-orbifolds.
This result has actually been part of the folklore for several years, but it turns
out that to state and prove it in detail one has to somewhat refine
the ideas underlying the proof in the manifold case.
This paper is devoted to these refinements.

We will state in the present
introduction two versions of the splitting
theorem that we will establish below, addressing the reader to
Section~\ref{statements:section} for the appropriate definitions.

\begin{teo}\label{basic:splitting:teo}
Let $X$ be a closed locally orientable $3$-orbifold. Suppose that $X$ does not contain
any bad $2$-suborbifold, and that every spherical $2$-suborbifold of $X$ is separating.
Then $X$ contains a finite system $\calS$ of spherical $2$-suborbifolds such that:
\begin{itemize}
\item No component of $X\setminus\calS$ is punctured discal.
\item If $Y$ is the orbifold obtained from $X\setminus\calS$ by capping the
boundary components with discal $3$-orbifolds, then the components of $Y$ are irreducible.
\end{itemize}
Any such system $\calS$ and the corresponding $Y$, if viewed as abstract collections
of $2$- and $3$-orbifolds respectively, depend on $X$ only.
\end{teo}

We inform the reader that this result has a dual interpretation in terms of connected sum,
discussed in detail in Section~\ref{statements:section}. We also note that the existence
part of Theorem~\ref{basic:splitting:teo} may seem to have a purely theoretical nature, which
is not pleasant in the context of $3$-dimensional topology, where a great emphasis
is typically given to algorithmic methods~\cite{Matveev:ATC3M}. The next longer statement emphasizes the
constructive aspects of the splitting.

\begin{teo}\label{splitting:intro:teo}
Let $X$ be a closed locally orientable $3$-orbifold. Suppose that $X$ does not contain
any bad $2$-suborbifold, and that every spherical $2$-suborbifold of $X$ is separating.
Starting with $\calS=\emptyset$ and $Y=X$, consider the (non-deterministic) process
described by the following steps:
\begin{itemize}
\item[1.] If all ordinary spherical $2$-suborbifolds of $Y$ are inessential, turn to Step 2.
Otherwise choose $\Sigma$ as one such $2$-suborbifold, redefine
$\calS$ as the given $\calS$ union $\Sigma$, redefine $Y$ as the orbifold obtained from
the given $Y$ by cutting along $\Sigma$ and capping the result by
discal $3$-orbifolds, and repeat Step 1.
\item[2.] If all cyclic spherical $2$-suborbifolds of $Y$ are inessential, turn to Step 3.
Otherwise choose $\Sigma$ as one such $2$-suborbifold, redefine $\calS$ and $Y$ as in Step 1,
and repeat Step 2.
\item[3.] If all vertex spherical $2$-suborbifolds of $Y$ are inessential, output $\calS$ and $Y$.
Otherwise choose $\Sigma$ as one such $2$-suborbifold, redefine $\calS$ and $Y$ as in Step 1,
and repeat Step 3.
\end{itemize}
Then the process is finite and the final $Y$ consists of irreducible $3$-orbifolds.
Moreover, the final $\calS$ and $Y$, viewed as abstract collections of $2$- and $3$-orbifolds
respectively, depend on $X$ only, not on the specific process leading to them.
In addition, at all steps of the process the essential spherical $2$-orbifold $\Sigma$ can be
chosen to be normal with respect to a triangulation of the given orbifold $Y$.
\end{teo}

For the reader familiar with the proof of the splitting theorem
for manifolds (see for instance~\cite{Hempel}) we will now
point out what goes wrong in the orbifold case, and sketch how we have
modified the argument. The main steps to show the existence of a
splitting of a manifold $M$ into prime ones are as follows:
\begin{itemize}
\item[(a)] Show that if $M$ contains a non-separating sphere then it splits as the connected sum of
$S^2\times S^1$ or $S^2\timtil S^1$ and some $M'$, and the first Betti number of $M'$ is smaller than that of $M$. Deduce
that one can restrict to $M$'s in which all spheres are separating.
\item[(b)] Define a family of spheres
$\calS\subset M$ to be \emph{essential} if
no component of $M\setminus \calS$ is a punctured disc.
Show that every essential $\calS$ can be replaced by one having the same number of components
and being normal with respect to any given triangulation of $M$. Deduce that all maximal
essential systems of spheres are finite.
\item[(c)] Prove that if $\calS$ is a maximal essential system then the manifolds obtained
by capping off the components of $M\setminus\calS$ are irreducible.
\end{itemize}
Uniqueness is more elaborate to prove, but its validity only depends on the
following fact, that we express using the dual viewpoint of connected sum:

\begin{itemize}
\item[(d)] The sphere $S^3$ is the only identity element of the operation of connected sum.
No connected sum with $S^3$ should be employed to realize a manifold.
\end{itemize}

Turning to orbifolds, we can show that none of the points (a)-(d) has a straight-forward
extension:
\begin{itemize}
\item[($\neg$a)] The existence of a non-separating spherical $2$-suborbifold does not imply the existence
of an essential separating one (see Fig.\ref{nonsep:fig} below).
In particular, there are infinitely many 3-orbifolds containing some non-separating
but no essential separating spherical $2$-orbifold (see Section~\ref{prime:section}).
\end{itemize}
The difficulty given by non-separating spherical $2$-orbifolds is mentioned in~\cite{BMP},
where the choice is made of cutting along these $2$-orbifolds too. The reasons of our choice
of excluding non-separating spherical $2$-orbifolds \emph{tout-court} will be discussed
in Section~\ref{existence:section}. We only want to mention here that the techniques of the present paper
do not seem to adapt easily to include the case of splitting along non-separating spherical $2$-orbifolds.

Turning to point (b), we follow~\cite{BMP} and define
a family $\calS$ of spherical $2$-suborbifolds of a 3-orbifold $X$ to be \emph{essential} if
no component of $X\setminus \calS$ is a punctured discal $3$-orbifold. We now have:

\begin{itemize}
\item[($\neg$b)] If $\calS$ is essential and $\calS',\calS''$ are obtained
from $\calS$ by one of the compression moves necessary to normalize $\calS$,
then neither of $\calS'$ and $\calS''$ need be essential (see Fig.~\ref{essnocompr:fig} below).
\end{itemize}

This fact could be viewed as a merely technical difficulty, showing that the theory
of normal surfaces is inadequate to prove finiteness of maximal essential systems.
The following fact is however more striking:

\begin{itemize}
\item[($\neg$c)] If $\calS$ is a maximal essential system of spherical $2$-suborbifolds of $X$,
and $Y$ is obtained by cutting $X$ along $\calS$ and capping the result, the components
of $Y$ need not be irreducible (see Fig.~\ref{essnoirred:fig} below).
\end{itemize}

Given this drawback of the notion of essential system, we have refined it
to a certain notion of \emph{efficient splitting system}, where we impose artificially
that the capped components of the complement should be irreducible. The proof of the existence of efficient
splitting systems then requires some more efforts than in the manifold case.
In addition, finiteness of the splitting process in Theorem~\ref{splitting:intro:teo}
does not follow directly from the existence of finite efficient
splitting systems.

Concerning uniqueness, we have the following:

\begin{itemize}
\item[($\neg$d)] There are three different types of connected sums of orbifolds, and each has its own
identity, different from the other ones. In a sequence of connected sums of irreducible $3$-orbifolds,
the property that a certain sum is ``trivial'' depends on the order in which the sums are performed.
\end{itemize}

The same notion of efficient system used to prove the existence of a splitting actually
allows to deal also with this problem, and hence to prove uniqueness in
Theorem~\ref{basic:splitting:teo}. The same argument employed for uniqueness also
proves finiteness of the splitting process in Theorem~\ref{splitting:intro:teo}.

\paragraph{Acknowledgements}
The author warmly thanks Gw\'ena\"el Massuyeau,
Ser\-gei Matveev, Luisa Paoluzzi, and Joan Porti for helpful
discussions.

\section{Efficient connected sums and splitting systems}\label{statements:section}
In this section we introduce the necessary terminology and we provide
alternative statements of Theorems~\ref{basic:splitting:teo} and~\ref{splitting:intro:teo}.

\paragraph{Local structure of orbifolds}
We will not cover here the general theory of orbifolds, referring the reader to the
milestone~\cite{thurston:notes}, to the excellent and very recent~\cite{BMP},
and to the comprehensive bibliography of the latter.
We just recall that an orbifold of dimension $n$ is a topological space with a singular smooth
structure, locally modelled on a quotient of $\matR^n$ under the action of a finite
group of diffeomorphisms. We will only need to refer to the cases $n=2$ and $n=3$, and we
will confine ourselves to orientation-preserving diffeomorphisms. In addition, all our
orbifolds will be compact.

Under these assumptions one can see that a 2-orbifold $\Sigma$
consists of a compact \emph{support} surface $|\Sigma|$ together with a finite
collection $S(\Sigma)$ of points in the interior of $|\Sigma|$, the \emph{cone} points,
each carrying a certain \emph{order} in $\{p\in\matN:\ p\geqslant2\}$.

Analogously, a $3$-orbifold $X$ is given by a compact \emph{support} 3-manifold $|X|$ together
with a \emph{singular set} $S(X)$. Here $S(X)$ is a finite collection of circles and unitrivalent
graphs tamely embedded in $|X|$, with the univalent vertices
given by the intersection with $\partial|X|$. Moreover each component of $S(X)$ minus the vertices
carries an \emph{order} in $\{p\in\matN:\ p\geqslant2\}$, with the restriction
that the three germs of edges incident to each vertex should
have orders $(2,2,p)$, for arbitrary $p$, or $(2,3,p)$, for $p\in\{3,4,5\}$.

\paragraph{Bad, spherical, and discal orbifolds}
An \emph{orbifold-covering} is a map between orbifolds locally modelled on a map of
the form $\matR^n/_\Delta\to\matR^n/_\Gamma$, naturally defined whenever
$\Delta<\Gamma<{\rm Diff}_+(\matR^n)$. An orbifold is called \emph{good} when
it is orbifold-covered by a manifold, and \emph{bad} when it is not good.
In the sequel we will need the following easy result:

\begin{lemma}\label{bad:lem}
The only bad closed $2$-orbifolds are $(S^2;p)$, the $2$-sphere with one cone point of
order $p$, and $(S^2;p,q)$, the $2$-sphere with cone points of orders $p\neq q$.
\end{lemma}

We now introduce some notation and terminology repeatedly
used below. We define $\ordtred$ to be $D^3$, the \emph{ordinary discal} $3$-orbifold,
$\cyctred(p)$ to be $D^3$ with singular set
a trivially embedded arc with arbitrary order $p$,
and $\vertred(p,q,r)$ to be $D^3$ with singular set
a trivially embedded
``Y-graph'' with edges of orders $p,q,r$.
We will call $\cyctred(p)$ and $\vertred(p,q,r)$ respectively
\emph{cyclic discal} and \emph{vertex discal} $3$-orbifolds, and we will employ
the shortened notation $\cyctred$ and $\vertred$ to denote
cyclic and vertex discal $3$-orbifolds with generic orders.

We also define the \emph{ordinary, cyclic}, and \emph{vertex spherical} $2$-orbifolds,
denoted respectively by $\orddues$, $\cycdues(p)$, and $\verdues(p,q,r)$,
as the 2-orbifolds bounding the corresponding discal 3-orbifolds
$\ordtred$, $\cyctred(p)$, and $\vertred(p,q,r)$.
We also define the \emph{ordinary, cyclic}, and \emph{vertex spherical} $3$-orbifolds,
denoted respectively by $\ordtres$, $\cyctres(p)$, and $\vertres(p,q,r)$,
as the 3-orbifolds obtained by mirroring the corresponding discal 3-orbifolds
$\ordtred$, $\cyctred(p)$, and $\vertred(p,q,r)$ in their boundary.
The spherical 2- and 3-orbifolds with generic orders will be denoted by $\matS^2_*$ and $\matS^3_*$.

\paragraph{2-suborbifolds and irreducible 3-orbifolds}
We say that a $2$-orbifold $\Sigma$ is a \emph{suborbifold} of a $3$-orbifold $X$ if $|\Sigma|$
is embedded in $|X|$ so that $|\Sigma|$ meets $S(X)$ transversely (in particular,
it does not meet the vertices), and $S(\Sigma)$ is given precisely by $|\Sigma|\cap S(X)$,
with matching orders.

A spherical $2$-suborbifold $\Sigma$ of a $3$-orbifold $X$ is called \emph{essential}
if it does not bound in $X$ a discal $3$-orbifold.
A $3$-orbifold is called \emph{irreducible} if it does not contain any
bad 2-suborbifold and every spherical 2-suborbifold is inessential
(in particular, it is separating).

If a $3$-orbifold $X$ is bounded by spherical $2$-orbifolds, we can canonically
associate to $X$ a closed orbifold $\widehat{X}$ by attaching the appropriate discal
$3$-orbifold to each component of $\partial X$. We say that $\widehat{X}$ is obtained
by \emph{capping} $X$.

Whenever a $3$-orbifold is not irreducible, we can select an essential
spherical $2$-suborbifold $\Sigma$, split $X$ along $\Sigma$, and cap the result. A na\"if
statement of the theorem we want to prove would be that \emph{the successive process of splitting
and capping comes to an end in finite time, and the result is a collection of
irreducible $3$-orbifolds, which depends on $X$ only}. With the appropriate assumptions and
choices of the essential spherical $2$-orbifolds this is eventually true (Theorem~\ref{splitting:intro:teo})
but, just as in the case of manifolds, the statement and proof are best understood by performing
a simultaneous splitting along a system of spherical $2$-orbifolds, rather than a successive
splitting along single ones.

\paragraph{Connected sum}
To motivate the statement of the main result given below, we reverse our point of view.
The remark here is that, if $X$ can be turned into a collection of irreducible $3$-orbifolds
by the ``split and cap'' strategy mentioned above, then $X$ can be reconstructed from the
same collection by ``puncture and glue'' operations. As in the case of
manifolds, each such operation will be called a \emph{connected sum}. But the world
of orbifolds is somewhat more complicated than that of manifolds, so we need to give the
definition with a little care.

Let $X_1$ and $X_2$ be 3-orbifolds, and pick points $x_i\in |X_i|$, such that one of the following
holds:
\begin{enumerate}
\item[(i)] Both $x_1$ and $x_2$ are non-singular points.
\item[(ii)] Both $x_1$ and $x_2$ are singular but not vertex points, and they belong
to singular edges of the same order $p$.
\item[(iii)] Both $x_1$ and $x_2$ are vertex singular points, and the triple $p,q,r$
of orders of the incident edges is the same for $x_1$ and $x_2$.
\end{enumerate}
We can then remove from $|X_i|$ a regular neighbourhood of $x_i$, and
glue together the resulting orbifolds by a homeomorphism
which matches the singular loci. The result is a
$3$-orbifold $X$ that we call a \emph{connected sum} of $X_1$ and $X_2$.
Depending on which of the conditions (i), (ii), or (iii) is met, we call
the operation (and its result), a connected sum of \emph{ordinary} type,
of \emph{cyclic} type \emph{(of order $p$)}, or of \emph{vertex} type
\emph{(of order $p,q,r$)}.

\paragraph{Trivial connected sums}
In the case of manifolds, there is only one type of operation of connected sum,
and its only identity element is the ordinary $3$-sphere. For $3$-orbifolds
we have the easy:

\begin{rem}\emph{
The unique identity element of the ordinary connected sum (respectively, of the order-$p$ cyclic
connected sum, and of the
the order-$(p,q,r)$ vertex connected sum)
is $\ordtres$ (respectively, $\cyctres(p)$, and $\vertres(p,q,r)$).}
\end{rem}

We can then define to be \emph{trivial} an operation of connected sum with the
corresponding identity element. However, when one considers a sequence of
operations of connected sums, one must be careful, because a sum
which appears to be non-trivial in the first place may actually turn out to
be trivial later on, as the following example, suggested by Joan Porti, shows.

\begin{example}\label{joan:example}
\emph{Let $K$ be any non-trivial knot in $S^3$, and let $p,q\geqslant 2$ be distinct
integers. Construct a $3$-orbifold $X$ by first taking the ordinary connected sum
of $\cyctres(p)$ and $\cyctres(q)$, and then taking the order-$p$ cyclic connected sum
of the result with $(S^3,K_p)$, the orbifold with support $S^3$ and singular set $K$ of
order $p$. Then both sums are non-trivial. However, if we reverse
the order, we can realize the same $X$ as the order-$p$ cyclic sum of
$\cyctres(p)$ and $(S^3,K_p)$, followed by the ordinary sum with $\cyctres(q)$, and
the first sum is now trivial.}
\end{example}

This phenomenon is of course a serious obstacle to proving that each
$3$-orbifold has a unique expression as a connected sum of irreducible ones.
To overcome this obstacle, we need to redefine the notion of trivial connected sum.

\paragraph{Graphs and efficient connected sums}
We consider in this paragraph successive connected sums $X_0\#\ldots\# X_n$,
stipulating with this notation that the first sum is between $X_1$ and $X_0$, the second
sum is between $X_2$ and $X_0\# X_1$, and so on.
Note first that to perform the $j$-th sum we can always arrange the
puncture made to $X_0\#\ldots\# X_{j-1}$,
to be disjoint from the spherical $2$-orbifolds along which
the punctured $X_0,\ldots,X_{j-1}$ have been glued.
This implies that $X_j$ is glued to precisely one of $X_0,\ldots,X_{j-1}$.
We can then construct a graph associated to the connected sum, with nodes
the $X_j$'s, and edges labelled by the types of sums performed. Note
that this graph is actually a tree.

We will say that a connected sum $X_0\#\ldots\# X_n$ is \emph{efficient} if
the $X_j$'s are irreducible and
in its graph there is no node $\ordtres$ with an incident edge of ordinary
type, there is no node $\cyctres$ with an incident edge of cyclic
type, and there is no node $\vertres$ with an incident edge of vertex type, see
Fig.~\ref{inefficientsums:fig}.
    \begin{figure}
    \begin{center}
    \input{inefficientsums.pstex_t}
    \nota{Inefficient connected sums.} \label{inefficientsums:fig}
    \end{center}
    \end{figure}

One should notice that the graph associated to $X_0\#\ldots\# X_n$ is
not quite unique, since to construct it one needs to isotope each puncture
away from the previous ones, as described above.  However, one easily sees
that two graphs of $X_0\#\ldots\# X_n$ differ by moves as in Fig.~\ref{manygraphs:fig},
    \begin{figure}
    \begin{center}
    \input{manygraphs.pstex_t}
    \nota{Moves relating two graphs of the same connected sum.} \label{manygraphs:fig}
    \end{center}
    \end{figure}
and one deduces that the notion of efficient connected sum is indeed well-defined.
The following is an alternative statement of our main result:

\begin{teo}\label{connected:sum:teo}
Let $X$ be a closed locally orientable $3$-orbifold. Suppose that $X$ does not contain
any bad $2$-suborbifold, and that every spherical $2$-suborbifold of $X$ is separating.
Then $X$ can be realized as an efficient connected sum of irreducible $3$-orbifolds.
Any two realizations involve the same irreducible summands and the same types of sums
(including orders).
\end{teo}

\paragraph{Efficient splitting systems}
As in the case of manifolds, to establish the existence
part of Theorem~\ref{connected:sum:teo} one proceeds
with the original strategy of ``split and cap'' first mentioned above,
but one has to split symultaneously along a (finite) system of disjoint spherical $2$-orbifolds.

To this end we give the following definitions. We say that a $3$-orbifold is \emph{punctured discal}
if it is obtained from some $\matD^{\,3}_*$ by removing a regular neighbourhood of a finite set.
A system $\calS$ of spherical $2$-suborbifolds of $X$ is called \emph{essential} if no component
of $X\setminus\calS$ is punctured discal, and it is called
\emph{coirreducible} if all the components of $(X\setminus\calS)\widehat{\ }$
are irreducible. We call \emph{efficient splitting system}
a finite system of spherical $2$-suborbifolds of $X$ which is both essential and coirreducible.

\begin{rem}\emph{To each system $\calS$ of separating spherical $2$-suborbifolds
there corresponds a realization of $X$ as a connected sum of the components of
$(X\setminus\calS)\widehat{\ }$. Moreover $\calS$ is efficient if and only if the connected sum is.}
\end{rem}

The result we will establish directly below, which is obviously equivalent to
Theorems~\ref{basic:splitting:teo} and~\ref{connected:sum:teo},
and will be used to deduce Theorem~\ref{splitting:intro:teo}, is the following:

\begin{teo}\label{efficient:splitting:teo}
Let $X$ be a closed locally orientable $3$-orbifold. Suppose that $X$ does not contain
any bad $2$-suborbifold, and that every spherical $2$-suborbifold of $X$ is separating.
Then $X$ admits efficient splitting systems. Any two splitting systems coincide
as abstract collections of spherical $2$-orbifolds, and the capped components of their
complements also coincide.
\end{teo}

\section{Comparison with the manifold case}\label{comparison:section}
We illustrate in greater detail in this section why it has been necessary to
substantially modify the proof of the splitting theorem in passing from
the manifold to the orbifold case.

Concerning the issue of non-separating spherical $2$-orbifolds (point (a) of the Introduction)
we begin by recalling what is the point for an orientable manifold $M$.
If $\Sigma\subset M$ is a non-separating sphere, one chooses a simple closed curve $\alpha$
meeting $\Sigma$ once and transversely and one takes the boundary $\Sigma'$ of
a regular neighbourhood of $\Sigma\cup\alpha$. Then $\Sigma'$ is again a sphere,
it is separating, and one of the capped components of $M\setminus\Sigma'$ is
$S^2\times S^1$ or $S^2\timtil S^1$, so $M$ has one such connected summand.
The very basis of this argument breaks down for orbifolds, because if $\Sigma$
is a spherical $2$-orbifold of cyclic or vertex type, then $\Sigma'$ is not
spherical at all (see Fig.~\ref{nonsep:fig}).
    \begin{figure}
    \begin{center}
\input{nonsep.pstex_t}
    \nota{A non-separating spherical $2$-orbifold may not determine a separating one.} \label{nonsep:fig}
    \end{center}
    \end{figure}
We will devote Section~\ref{prime:section} to showing
that indeed there is a wealth of $3$-orbifolds which contain
non-separating spherical $2$-orbifolds but cannot be expressed as a non-trivial connected sum.

Turning to points (b) and (c) of the Introduction, recall that
a family $\calS$ of spherical $2$-suborbifolds of a 3-orbifold $X$ is called \emph{essential} if
no component of $X\setminus \calS$ is a punctured discal $3$-orbifold.
To prove finiteness of a maximal essential family using the analogue of orbifolds
of Haken's normal surface theory, one should be able to prove that the property of
being essential is stable under the ``normalization moves.'' Now one can see (as we will
in Section~\ref{uniqueness:section}), that the normalization moves boil down to isotopy
and compression along ordinary discs (\emph{i.e.} discs without singular points).
And we have the following:

\begin{example}\label{non:b:exa}
\emph{Let a $3$-orbifold $X$ contain a system $\calS$ of three cyclic spherical $2$-suborbifolds or orders
$p$, $p$, and $q$, for some $p,q\geqslant 2$. Suppose that one of the components of
$X\setminus\calS$ is as shown in Fig.~\ref{essnocompr:fig}
    \begin{figure}
    \begin{center}
    \input{essnocompr.pstex_t}
    \nota{Compression along an ordinary disc does not preserve essentiality.} \label{essnocompr:fig}
    \end{center}
    \end{figure}
(an ordinary connected sum of $\cyctres(p)$ and $\cyctres(q)$, with two cyclic punctures of order $p$ and
one of order $q$). Suppose that the other components of $X\setminus\calS$ are once-punctured
irreducible non-spherical $3$-orbifolds. Let $\calS'$ and $\calS''$ be the systems obtained from
$\calS$ by compression along the ordinary disc $D$ also shown in Fig.~\ref{essnocompr:fig}.}
Then $X$ does not contain any bad $2$-suborbifold, every spherical $2$-suborbifold of $X$ is separating,
$\calS$ is essential, $\calS'$ and $\calS''$ are systems of spherical $2$-orbifolds, and they are both
inessential in $X$.
\end{example}

The failure of the direct orbifold analogue of normal surfaces to establish existence of a maximal
essential system may appear to be a technical difficulty, only calling for a different proof
of existence. But the next example, which refers to point (c) in the Introduction, shows that the notion
itself of \emph{essential} system does not do the job it was designed for:

\begin{example}\label{non:c:exa}
\emph{Let a $3$-orbifold $X$ contain a system $\calS$ of one cyclic and two vertex spherical
$2$-suborbifolds or orders $p$ and $p,q,r$ (twice) for some admissible
$p,q,r\geqslant 2$. Suppose that one of the components of
$X\setminus\calS$ is as shown in Fig.~\ref{essnoirred:fig},
    \begin{figure}
    \begin{center}
    \input{essnoirred.pstex_t}
    \nota{The capped components of
    the complement of a maximal essential system need not be irreducible.} \label{essnoirred:fig}
    \end{center}
    \end{figure}
that is, $\vertres(p,q,r)$ with two vertex and one order-$p$ cyclic punctures).
Let $\calS'$ be given by the two vertex components of $\calS$.
Suppose that the other components of $X\setminus\calS$ are once-punctured
irreducible non-spherical $3$-orbifolds.
Then} $X$ does not contain any bad $2$-suborbifold, every spherical $2$-suborbifold of $X$ is separating,
$\calS'$ is a maximal essential system spherical $2$-orbifolds in $X$, but not all the capped
components of $X\setminus\calS'$ are irreducible.
\emph{To understand the last two assertions,
note that $\calS$ is not essential, because the $3$-orbifold in
Fig.~\ref{essnoirred:fig} is indeed punctured discal, but the cyclic
component of $\calS$ is essential after cutting along $\calS'$ and capping,
because it corresponds to a} cyclic \emph{connected sum with $\vertres(p,q,r)$, which
is non-trivial}.
\end{example}

It is perhaps worth pointing out the basic fact underlying the previous example, namely
that \emph{a punctured spherical $3$-orbifold is a discal one if and only if the
puncture has the same type as the original $3$-orbifold}.

Concerning uniqueness (point (d) of the Introduction), all we had to say
about its statement was said in Section~\ref{statements:section}. About its proof,
we only note that Milnor's one for manifolds was based on the same method proving
existence, \emph{i.e.} on the
notion of essential spherical system, so of course our proof will be somewhat different.

\section{Existence of efficient splitting systems and\\ algorithmic aspects}\label{existence:section}
In this section we prove the existence parts of
Theorems~\ref{efficient:splitting:teo} and~\ref{connected:sum:teo}.
We do so by showing that the process described
in Theorem~\ref{splitting:intro:teo} can be carried out with all the components of the
splitting system simultaneously normal with respect to a  fixed triangulation.

\paragraph{Triangulations and normal $2$-orbifolds}
We slightly relax the definition given in~\cite{BMP}, and define a \emph{triangulation}
of a $3$-orbifold $X$ as a triangulation of $|X|$ which contains $S(X)$ as a subcomplex of its
one-skeleton. We then define a $2$-suborbifold of $X$ to be normal with respect to a triangulation
if it intersects each tetrahedron in a union of squares and triangles, just as in the case
of a surface in a triangulated manifold. The following fundamental
result holds:

\begin{prop}\label{normal:finiteness:prop}
Let $X$ be a triangulated $3$-orbifold. Then there exists a number $m$
such that, if $\calF$ is a normal surface with respect to the given triangulation
and $\calF$ has more than $m$ connected components, then
at least one of the components of $X\setminus\calF$ is a product orbifold.
\end{prop}

This result is established exactly as in the case of manifolds, so we omit the proof.
One only needs to notice that if a region of $X\setminus\calF$ is a product as a manifold
then it is also a product as an orbifold.

\begin{figure}
    \begin{center}
    \mettifig{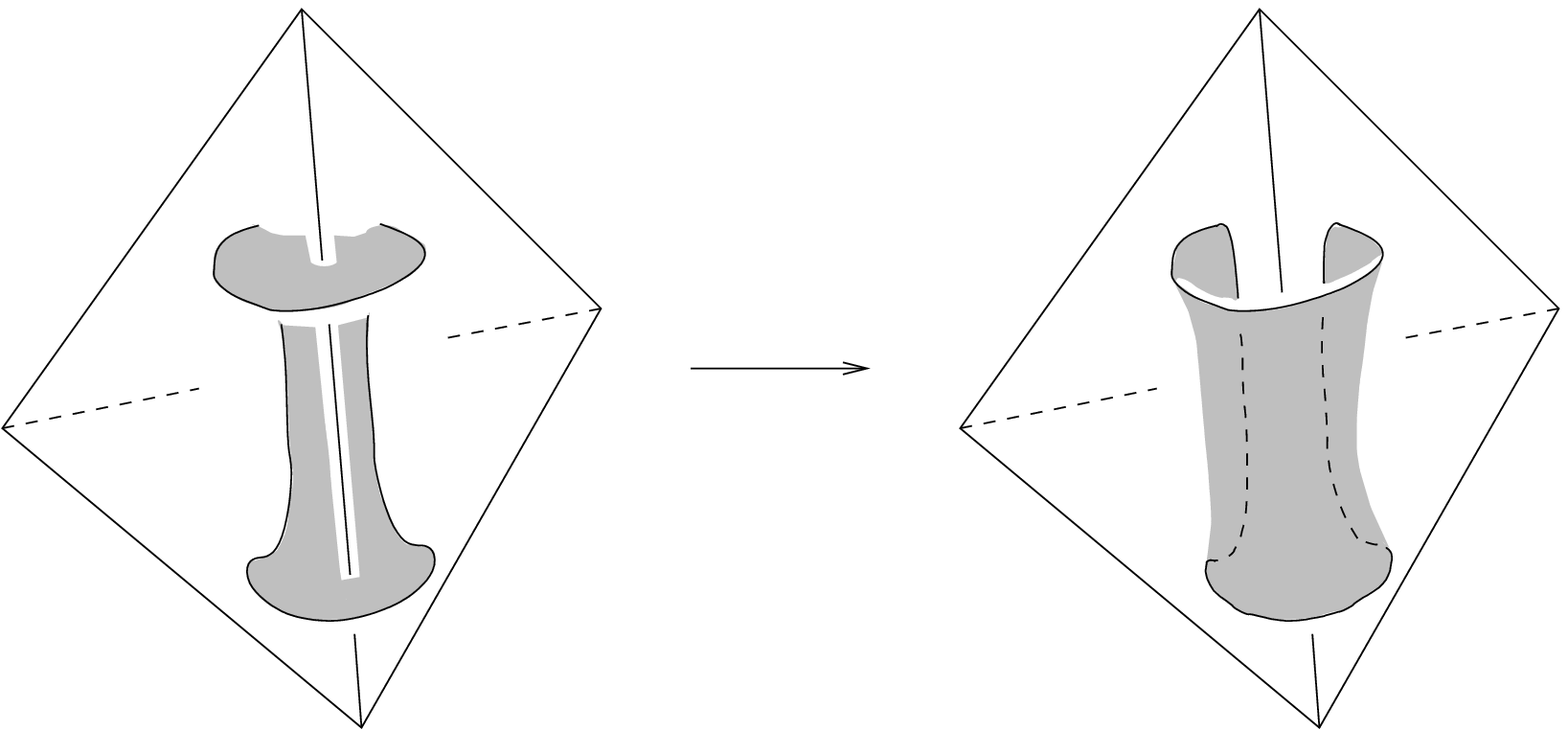,width=9cm}
    \nota{Third normalization move.} \label{secondnorm:fig}
    \end{center}
    \end{figure}

\paragraph{Normalization moves}
We now recall that, in a manifold with a fixed triangulation $T$, there is a general strategy to replace a surface
by a normal one, as explained for instance in~\cite{Matveev:ATC3M}. When the surface is a sphere $\Sigma$,
this strategy consists (besides isotopy relative to $T$) of three ``normalization''
moves. \emph{The first (respectively, second) move consists in compressing $\Sigma$ along a disc $D$ contained
in the interior of a tetrahedron (respectively, triangle) of $T$,
thus getting two spheres $\Sigma'$ and $\Sigma''$, and choosing either $\Sigma'$ or $\Sigma''$.}
The third move is used to remove the double intersections of the normal discs with the edges, and
it is illustrated in Fig.~\ref{secondnorm:fig}.
As a matter of fact, this move can also be seen as compression along a disc
followed by the choice of one of the resulting spheres. More precisely, one
compresses along the disc in Fig.~\ref{secondnorm:fig}-right
and one discards the sphere given by the union of the two discs in Fig.~\ref{secondnorm:fig}.
Therefore, \emph{the third move consists of the compression
along a disc disjoint from the $1$-skeleton of $T$, followed by the elimination
of a sphere which is the boundary of a regular neighbourhood of an interior point of an
edge of $T$.}

\paragraph{Stepwise normalization}
As already mentioned, to prove the existence of an efficient splitting system
we will construct a system of spherical $2$-orbifolds along the steps
of Theorem~\ref{splitting:intro:teo}, making sure each new component
of the system is normal with respect to a fixed triangulation of
the original orbifold. To show that this is possible, we start with
an easy lemma, that we will use repeatedly below.
For its statement, we stipulate that a spherical 2-orbifold of cyclic type is
\emph{more complicated} than one of ordinary type, and a
that a spherical 2-orbifold of vertex type is
\emph{more complicated} than one of ordinary or cyclic type.

\begin{lemma}\label{discal:cap:lem}
Let $\Delta$ be a punctured discal $3$-orbifold, and let $\Sigma$ be a component
of $\partial\Delta$. Suppose that no component of $\partial\Delta$ is more
complicated than $\Sigma$. Then capping all the components of $\partial\Delta$
except $\Sigma$ we get a discal $3$-orbifold.
\end{lemma}

The next two results are the core of our argument on stepwise normalization.
The first one will be used again in the next section. From now on
we fix a $3$-orbifold $X$ and we suppose that $X$ does not contain
any bad $2$-suborbifold, and that every spherical $2$-suborbifold of $X$ is separating.

\begin{prop}\label{ess:ess:prop}
Let $\calS$ be an essential system of spherical $2$-suborbifolds of $X$.
Let $Y$ be a component of $X\setminus\calS$ and $\Sigma$ be a spherical $2$-suborbifold of $Y$.
Suppose that:
\begin{itemize}
\item $\Sigma$ is essential in $\hatY$.
\item All the components of $\calS$ are at most as complicated as $\Sigma$.
\end{itemize}
Then $\calS\cup\Sigma$ is essential.
\end{prop}

\begin{proof}
Suppose that some component $Z$ of $X\setminus(\Sigma\cup\calS)$ is punctured discal.
The assumptions that $\calS$ is essential and that $\Sigma$ is separating imply that
$Z$ is contained in $Y$ and $\Sigma$ is a boundary component of $Z$.
By Lemma~\ref{discal:cap:lem} and the second assumption,
capping all the components of $\partial Z$ but $\Sigma$ we get a discal $3$-orbifold, which
implies that $\Sigma$ bounds a discal $3$-orbifold in $\widehat{Y}$: a contradiction to the first assumption.
\end{proof}

\begin{prop}\label{step:norm:prop}
Fix a triangulation $T$ of $X$.
Let $\calS$ be an essential system of spherical $2$-suborbifolds
in normal position with respect to $T$.
Let $Y$ be a component of $X\setminus\calS$ and $\Sigma$ be a spherical $2$-suborbifold of $Y$.
Suppose that:
\begin{itemize}
\item $\Sigma$ is essential in $\hatY$.
\item All the spherical $2$-suborbifolds of $\hatY$
strictly less complicated than $\Sigma$ are inessential in $\hatY$.
\end{itemize}
Then $\Sigma$ can be replaced by a spherical $2$-orbifold satisfying the same
properties and such that $\calS\cup\Sigma$ is in normal position with respect to $T$.
\end{prop}

\begin{proof}
We must show that a normalization move applied to $\Sigma$ can be realized without touching $\calS$,
without changing the type of $\Sigma$, and
preserving the property that $\Sigma$ is essential $\hatY$.
The first assertion is a straight-forward consequence of the description of the normalization moves
given above. For the second and third assertions, we first note that in all cases the compression of a normalization
move takes place along a disc disjoint from the $1$-skeleton of $T$, whence along an ordinary (non-singular)
disc $D$ which meets $\Sigma\cup\calS$ in $\partial D\subset \Sigma$ only.
To fix the notation, let $D'$ and $D''$ be the discs bounded by $\partial D$ on $\Sigma$, and define
$\Sigma'=D\cup D'$ and $\Sigma''=D\cup D''$.
By the assumptions that $X$ contains no bad 2-suborbifolds and that $D$ is an ordinary disc,
we can (and will) assume up to permutation that $\Sigma''$ is an ordinary sphere,

We now claim that $\Sigma'$ or $\Sigma''$ is essential
in $\hatY$. Before proving the claim, let us see how it implies the conclusion.
Since $\Sigma''$ is ordinary, we have $\Sigma'\cong\Sigma$.
If $\Sigma'$ is essential, we are done. Otherwise $\Sigma''$ is essential (and ordinary),
so $\Sigma$ is also ordinary by the second assumption, and we can switch $\Sigma'$ and $\Sigma''$.

To prove the claim,
suppose by contradiction that $\Sigma'$ and $\Sigma''$ are both inessential in $\hatY$.
Since $\Sigma''$ is ordinary, we deduce that $\Sigma$ bounds in $\widehat{Y}$
either a discal $3$-orbifold glued along a 2-disc to an ordinary 3-disc, or a 3-orbifold
obtained from a discal one by removing an ordinary 3-disc which intersects the boundary in a 2-disc.
The result is in both cases a discal $3$-orbifold (isomorphic to the original one),
whence the conclusion.
\end{proof}

The result just proved easily implies the following:

\begin{cor}
If $X$ is triangulated then one can carry out the process of Theorem~\ref{splitting:intro:teo}
with the splitting system $\calS$ in normal position at each step.
\end{cor}

This corollary, together with
Propositions~\ref{normal:finiteness:prop} and~\ref{ess:ess:prop}
immediately yields the result we were heading for:

\begin{cor}
$X$ admits finite efficient splitting systems.
\end{cor}

As promised in the Introduction, we explain here why it is crucial for us to assume that
all the spherical $2$-suborbifolds of $X$ are separating. First of all, the definition itself
of essential system would not be appropriate without this assumption, because a
non-separating spherical 2-orbifold could constitute an inessential system.
In addition, the assumption was used in a fundamental way in the proof of Proposition~\ref{ess:ess:prop}.

Proposition~\ref{step:norm:prop} also has the next consequence, which proves
the last assertion of Theorem~\ref{splitting:intro:teo}:

\begin{cor}
Let $Z$ be a triangulated orbifold and let $\Sigma$ be an essential spherical
$2$-suborbifold of $Z$. Suppose that all the spherical $2$-suborbifolds
of $Z$ strictly less complicated than $\Sigma$ are inessential.
Then there exists an essential spherical
$2$-suborbifold of $Z$ of the same type as $\Sigma$ and in normal position
with respect to the fixed triangulation.
\end{cor}

\section{Uniqueness of the splitting}\label{uniqueness:section}
In this section we conclude the proofs of
Theorems~\ref{efficient:splitting:teo},~\ref{connected:sum:teo}, and~\ref{splitting:intro:teo}.
As above we fix a $3$-orbifold $X$ and we suppose that $X$ does not contain
any bad $2$-suborbifold, and that every spherical $2$-suborbifold of $X$ is separating.

\paragraph{Gluing of discal $3$-orbifolds}
Our uniqueness result uses the following two easy technical lemmas.
We include a proof of the second one, leaving to the reader the even easier
proof of the first one.

\begin{lemma}\label{disc:attach:lem}
Let $W$ be a $3$-orbifold and let $D\subset\partial|W|$ be a disc. Let $\Delta$ be a punctured
discal $3$-orbifold. Let $\Sigma$ be a component of $\partial\Delta$ of the same type as $\Delta$.
Split $|\Sigma|$ along a loop into two discs $D'$ and $D''$, and suppose that $D,D',D''$
are isomorphic to each other, with at most one singular point. Let $W'$ be obtained
from $W$ by gluing $\Delta$ along a homeomorphism $D'\to D$. Then:
\begin{enumerate}
\item $W'$ is punctured discal if and only if $W$ is.
\item $\widehat{W}'=\widehat{W}$.
\end{enumerate}
\end{lemma}

\begin{lemma}\label{sphere:attach:lem}
Let $W$ be a $3$-orbifold with a spherical
boundary component $\Sigma$. Let $\Delta$ be a punctured
discal $3$-orbifold with at least two boundary components $\Sigma',\Sigma''$
isomorphic to $\Sigma$, and all other components not more complicated than $\Sigma'$ and $\Sigma''$.
Let $W'$ be obtained
from $W$ by gluing $\Delta$ along a homeomorphism $\Sigma'\to\Sigma$. Then:
\begin{enumerate}
\item $W'$ is punctured discal if and only if $W$ is.
\item $\widehat{W}'=\widehat{W}$.
\end{enumerate}
\end{lemma}

\begin{proof}
We begin by (2). Since $\Sigma'$ and $\Sigma''$ are at least as complicated as all the
other components of $\partial\Delta$,
capping the components different from $\Sigma'$ and $\Sigma''$
we get a product orbifold, which implies the conclusion at once.

Let us turn to (1). To be discal for a $3$-orbifold $Z$ means that $\widehat{Z}$ is some
$\matS^3_{\rm t}$ and some component of $\partial Z$, or more precisely the most complicated one,
is $\matS^2_{\rm t}$. Now we already know that $\widehat{W}'=\widehat{W}$,
and the assumptions imply that the most complicated component of
$\partial{W}'$ has the same type as the most complicated component of
$\partial{W}$, whence the conclusion.
\end{proof}

\paragraph{Finiteness and uniqueness: scheme of the proof}
To conclude the proofs of the results stated so far we must establish
the finiteness of the splitting process of
Theorem~\ref{splitting:intro:teo}
and the uniqueness part of Theorems~\ref{splitting:intro:teo},~\ref{connected:sum:teo},
and~\ref{efficient:splitting:teo}. The former easily follows from the
next result, to be proved below, and Proposition~\ref{ess:ess:prop}.

\begin{prop}\label{lexi:prop}
Let $\calS$ and $\calS'$ be finite systems of spherical $2$-suborbifolds of $X$. Denote
by $n_{\rm t}$ (respectively, $n'_{\rm t}$) the number of components of $\calS$
(respectively, $\calS'$) of type ${\rm t}$. Suppose that $\calS$ is coirreducible
and $\calS'$ is essential. Then $(n'_{\rm o},n'_{\rm c},n'_{\rm v})\leqslant
(n_{\rm o},n_{\rm c},n_{\rm v})$ in lexicographical order.
\end{prop}

Uniqueness easily follows from the next result (for Theorem~\ref{splitting:intro:teo}
one only has to note that the final $\calS$ is obviously coirreducible,
and essential by Proposition~\ref{ess:ess:prop}):

\begin{prop}\label{uniqueness:prop}
Let $\calS$ and $\calS'$ be efficient spherical splitting systems for $X$.
Then $\calS\cong\calS'$ and $(X\setminus\calS)\widehat{\ }\cong(X\setminus\calS')\widehat{\ }$
as abstract collections of $2$- and $3$-orbifolds respectively.
\end{prop}

\paragraph{Moves for spherical systems}
We will prove Propositions~\ref{lexi:prop} and~\ref{uniqueness:prop}
in a unified manner. The underlying idea in both cases is that by certain moves we can modify $\calS$ into a superset
of $\calS'$, leaving unaffected the properties of $\calS$ and the topological types of
$\calS$ and $(X\setminus\calS)\widehat{\ }$ as abstract orbifolds. There are two types of move
we will need, that we now describe. These moves apply to any spherical system $\calS\subset X$, under
the appropriate assumptions.

\emph{Move $\alpha$}. Suppose we have a component $\Sigma_0$ of $\calS$ and a disc $D_1$ properly embedded
in a component $Y_0$ of $X\setminus\calS$, with $\partial D_1\subset\Sigma_0$.
Let $\partial D_1$ split $\Sigma_0$ into discs $D_0$ and $D$, and suppose that $D_0$ and $D_1$
both have at most one singular point. Assume that the component containing $D_0$
of $Y_0\setminus D_1$ is a punctured disc $\Delta$, and that $\partial\Delta$ does not
have components more complicated than $D_0\cup D_1$. Then we replace $\Sigma_0=D\cup D_0$ by
$\Sigma_1=D\cup D_1$. See Fig.~\ref{alphamove:fig}.
    \begin{figure}
    \begin{center}
    \input{alphamove.pstex_t}
    \nota{The two instances where move $\alpha$ applies, with $\widehat\Delta=\ordtres$ (left)
    and $\widehat\Delta=\cyctres$ (right).} \label{alphamove:fig}
    \end{center}
    \end{figure}

\emph{Move $\beta$}. Suppose we have a component $\Sigma_0$ of $\calS$ and a
spherical $2$-orbifold $\Sigma_1$ contained
in a component $Y_0$ of $X\setminus\calS$, with $\Sigma_0\subset\partial Y_0$.
Suppose that the component containing $\Sigma_0$ of $Y_0\setminus\Sigma_1$
is a punctured disc $\Delta$, that $\Sigma_0$ and $\Sigma_1$ have the same type, and
that no other component of $\partial\Delta$ is more complicated than $\Sigma_0$ and $\Sigma_1$. Then we replace
$\Sigma_0$ by $\Sigma_1$. See Fig.~\ref{betamove:fig}.
    \begin{figure}
    \begin{center}
    \input{betamove.pstex_t}
    \nota{The three instances where move $\beta$ applies, with $\widehat\Delta=\ordtres$ (left),
    $\widehat\Delta=\cyctres$ (centre), and $\widehat\Delta=\vertres$ (right).} \label{betamove:fig}
    \end{center}
    \end{figure}

\begin{lemma}\label{abmoves:lem}
Under any move of type $\alpha$ or $\beta$ the following happens:
\begin{itemize}
\item The type of $\calS$ (including orders) does not change.
\item The type of $(X\setminus\calS)\widehat{\ }$ does not change (in particular, coirreducibility
is preserved).
\item Essentiality is preserved.
\end{itemize}
\end{lemma}

\begin{proof}
The first assertion is proved by direct examination of Figg.~\ref{alphamove:fig}
and~\ref{betamove:fig}. To prove the second and third assertions we
note the only components of $X\setminus\calS$ modified by the move are $Y_0$ and the other
component $Z_0$ incident to $\Sigma_0$, and these components are replaced by the components
$Y_1$ and $Z_1$ incident to $\Sigma_1$. The relations between them are described as follows:
\begin{eqnarray*}
{\rm move\ }\alpha:  &  Y_0=Y_1\cup_{D_1}\Delta, & Z_1=Z_0\cup_{D_0}\Delta; \\
{\rm move\ }\beta:  & Y_0=Y_1\cup_{\Sigma_1}\Delta, & Z_1=Z_0\cup_{\Sigma_0}\Delta.
\end{eqnarray*}
It is now easy to check that in all cases we are under the assumptions of Lemmas~\ref{disc:attach:lem}
and~\ref{sphere:attach:lem}, whence the conclusion at once.\end{proof}

\paragraph{Finiteness and uniqueness: conclusion}
As already mentioned, to prove Propositions~\ref{lexi:prop} and~\ref{uniqueness:prop}
we want to replace $\calS$ by a superset of $\calS'$. The first step is as follows:

\begin{prop}\label{no:intersection:prop}
Suppose that $\calS$ is coirreducible and $\calS'$ is essential.
Then, up to applying moves $\alpha$ and isotopy to $\calS$, we can suppose that
$\calS\cap\calS'=\emptyset$.
\end{prop}

\begin{proof}
Let us first isotope $\calS$ so that it is transversal to $\calS'$.
Of course it is sufficient to show that if $\calS\cap\calS'$ is non-empty then
the number of its components can be reduced by a move $\alpha$ and isotopy.

Suppose some component $\Sigma'$ of $\calS'$ meets $\calS$, and note
that the finite set of loops $\calS\cap\Sigma'$ bounds at least two innermost discs in $\Sigma'$.
Since $\Sigma'$ has at most three singular points, we can choose an innermost disc $D_1$
which contains at most one singular point. Let $Y_0$ be the component of
$X\setminus\calS$ which contains $D_1$, let $\Sigma_0$ be the component of $\calS$
containing $\partial D_1$, and denote by $D_0$ and $D$ the discs into which
$\partial D_1$ cuts $\Sigma_0$, with $D_0$ at most as complicated as $D$.

The assumption that $X$ contains no bad 2-suborbifolds
easily implies that $D_0\cup D_1$ is a spherical $2$-orbifold of ordinary or cyclic type.
Coirreducibility of $\calS$ then implies that
$D_0\cup D_1$ bounds an ordinary or cyclic discal $3$-orbifold $B$ in $\hatY_0$.
Now we either have $B\cap D=\emptyset$ or $B\supset D$. In the first case
it is easy to see that we are precisely in the situation where we can apply a move $\alpha$,
after which a small isotopy allows to reduce the number of components of $\calS\cap\calS'$.
In the second case $D_1\cup D$ bounds in $\hatY_0$ one of the manifolds obtained from $B$
by cutting along $D$. Since $B$ is an ordinary or cyclic disc we deduce that $D$ is isomorphic
to $D_1$ and that $D\cup D_1$ bounds a ball in $\hatY_0$ disjoint from $D_0$. We can then interchange the roles
of $D$ and $D_0$ and apply $\alpha$ also in this case, whence the conclusion.
\end{proof}

Propositions~\ref{lexi:prop} and~\ref{uniqueness:prop} now easily follow from
Proposition~\ref{no:intersection:prop} and the next result:

\begin{prop}\label{ord:fix:prop}
Suppose that $\calS$ is coirreducible and $\calS'$ is essential.
Assume that each component of $\calS'$ is either contained in $\calS$
or disjoint from $\calS$.
Then:
\begin{itemize}
\item Up to applying moves $\beta$ and isotopy to $\calS$, we can suppose that
each ordinary component of $\calS'$ is contained in $\calS$.
\item If $\calS$ and $\calS'$ have the same ordinary components,
up to applying moves $\beta$ and isotopy to $\calS$, we can suppose that
each cyclic component of $\calS'$ is contained in $\calS$.
\item If $\calS$ and $\calS'$ have the same ordinary and cyclic components,
up to applying moves $\beta$ and isotopy to $\calS$, we can suppose that
each vertex component of $\calS'$ is contained in $\calS$.
\end{itemize}
\end{prop}

\begin{proof}
In all the cases, we suppose that a certain component $\Sigma_1$ of $\calS'$ of
the appropriate type is not contained $\calS$, and we show that a move $\beta$
can be applied to $\calS$ which replaces some component $\Sigma_0$ by $\Sigma_1$.
The conclusion then easily follows by iteration. In all cases we
denote by $Y_0$ the component of $X\setminus\calS$ which contains $\Sigma_1$,
and we choose $\Sigma_0$ among the components of $\partial Y_0$.

We start with the first assertion, so we assume $\Sigma_1$ is ordinary.
Coirreducibility of $\calS$ implies that $\Sigma_1$ bounds an ordinary 3-disc $B$ in
$\hatY_0$. Let us choose a component of $\calS'$ contained in $B$, not contained in
$\calS$, and innermost in $B$ with respect to these properties, and let
us rename it $\Sigma_1$. (Note that the new $\Sigma_1$ may not be parallel to the old one in $X$,
even if it is in $B$). Of course $\Sigma_1$ is again an ordinary sphere, and it bounds
in $Y$ a punctured ordinary disc $\Delta$ such that all the components of
$\partial\Delta$ except $\Sigma_1$ belong to $\calS$. Essentiality of $\calS'$ now implies that
there must be a component $\Sigma_0$ of $\partial\Delta$ which does not belong to $\calS'$.
We have created a situation where a move $\beta$ can be applied, so we are done.

The proof is basically the same for the second and third assertions.
We only need to note that $B$ always has the same type as $\Sigma_1$ and that
the component $\Sigma_0$ of $\partial\Delta$ not belonging to $\calS'$
cannot be of type strictly less complicated than $\Sigma_1$, by the assumption that $\calS$ and
$\calS'$ have the same components of type strictly less complicated than that of $\Sigma_1$.
\end{proof}

\section{Prime non-irreducible 3-orbifolds}\label{prime:section}

In this section we treat point (d) of the Introduction, showing that there
are infinitely many $3$-orbifolds with non-separating but no essential separating
spherical $2$-suborbifolds. Our examples arose from discussions with Luisa Paoluzzi.

As in the case of manifolds, one can define a prime orbifold as one that cannot
be expressed as a non-trivial connected sum, \emph{i.e} an orbifold in which every
\emph{separating} spherical $2$-suborbifold bounds a discal $3$-suborbifold.
Of course every irreducible orbifold is also prime, but we
will show now that there are infinitely many primes which
are not irreducible. This appears as a sharp difference between the orbifold
and the manifold case (in which the only non-irreducible primes are the two $S^2$-bundles over $S^1$).
This difference represents a serious obstacle to promoting
Theorem~\ref{connected:sum:teo} to a unique decomposition theorem for orbifolds into
prime ones.

\begin{prop}\label{non:irred:prop}
Let $X$ be a $3$-orbifold bounded by a sphere with $4$ cone points of the same
order $p$.  Assume that:
\begin{itemize}
\item $X$ is irreducible.
\item $X$ is not the $3$-ball with two parallel unknotted singular arcs.
\item In $X$ there is no sphere with three singular points which meets
$\partial X$ in a disc with two singular points.
\end{itemize}
Consider the orbifold with support $S^2\times S^1$ and singular set given by two circles
of order $p$ parallel to the  factor $S^1$. Remove from it
a regular neighbourhood of an arc
which joins the singular components
within a level sphere $S^2\times\{*\}$, and call $Y$ the result.
Let $Z$ be obtained by attaching $X$ and $Y$ along their boundary spheres, matching
the cone points. Then $Z$ is prime but not irreducible.
\end{prop}

\begin{proof}
A level sphere $S^2\times\{*\}$ contained in $Y$ gives a non-separating spherical $2$-suborbifold
of $Z$, so $Z$ is not irreducible.

Let us consider a separating spherical $2$-orbifold $\Sigma$ in $Z$, and the
intersection of $\Sigma$ with the sphere $\Sigma_0$
along which $X$ and $Y$ have been glued together. If $\Sigma\cap \Sigma_0$ is empty then
$Y$ bounds a spherical $3$-orbifold, because $X$ and $Y$ are irreducible.
Let us then assume $\Sigma\cap \Sigma_0$ to be transverse and minimal up to isotopy of $\Sigma$.
Considering the pattern of circles $\Sigma\cap \Sigma_0$ on $\Sigma$ we see that $\Sigma$ contains
at least two innermost discs, which are therefore contained either in $X$ or in $Y$.
Since $\Sigma$ has at most three singular points, we can find one such disc $D$ with
$i\leqslant 1$ points. Now consider the loop $\partial D$ on $\Sigma_0$, and
recall that $\Sigma_0$ has $4$ singular points. Then
$\partial D$ encircles $e\leqslant 2$ singular points. We will now show that
both if $D\subset X$ and if $D\subset Y$ the conditions $i\leqslant 1$ and $e\leqslant 2$
are impossible, whence the conclusion.

Suppose first that $D\subset X$. Using the fact that $X$ is irreducible, the
cases $e=i=0$ and $e=i=1$ both contradict the minimality of $\Sigma\cap \Sigma_0$. The cases $e=0,\ i=1$
and $e=1,\ i=0$ both lead to a bad 2-suborbifold of $X$, which does not exist.
In the case $e=2,\ i=0$, using irreducibility twice, we see that $X$ must be
the $3$-ball with two parallel unknotted singular arcs, but this orbifold was specifically forbidden,
so indeed we get a contradiction.
The case $e=2,\ i=1$ was also specifically forbidden.

Suppose now that $D\subset Y$. Recalling that $D$ is contained in a sphere which
is separating in $Z$, we see that $D$ cobounds a ball with a disc contained in $\partial Y$.
Using this fact we can prove again that all the possibilities with
$e\leqslant 2$ and $i\leqslant 1$ are impossible. As above,
$e=i=0$ and $e=i=1$ contradict minimality. The case $e=2,\ i=0$ would imply that one of the
singular arcs of $Y$ can be homotoped to $\partial Y$, which is not the case. The other
cases are absurd for similar reasons.
\end{proof}

\paragraph{Examples}
Whenever $k\geqslant1$ and the cone angles are chosen in an admissible fashion,
the ``stairway'' of Fig.~\ref{stairs:fig} contained in the 3-ball
    \begin{figure}
    \begin{center}
    \input{stairs.pstex_t}
    \nota{An example or orbifold $X$ such that $Z=X\cup Y$ is prime but not irreducible.} \label{stairs:fig}
    \end{center}
    \end{figure}
provides an example of orbifold $X$ satisfying the assumptions of Proposition~\ref{non:irred:prop}.
Note that for such an $X$ the corresponding prime non-irreducible orbifold $Z$ is always supported
by $S^2\times S^1$, so $|Z|$ is prime as a manifold. But we can also construct
examples of $X$ to which Proposition~\ref{non:irred:prop}
applies, leading to an orbifold $Z$ with non-prime $|Z|$.
For instance, one can take an orbifold constructed exactly as $Y$ was constructed
from $S^2\times S^1$, starting instead with $F\times S^1$, where $F$ is
a closed orientable surface $F$ of positive genus.

\section{3-orbifolds with ordinary and cyclic splitting}\label{ord:cyc:section}
Our interest in the spherical splitting of $3$-orbifolds arose from our program of
extending Matveev's complexity theory~\cite{Bruno:survey,Matveev:AAM} from
the manifold to the orbifold context. We actually succeeded~\cite{orb:complexity} in generalizing
the definition and some of the most significant results,
investigating in particular the behaviour under connected sum. However it turns out
that one does not have in general plain additivity as in the manifold case. We can
only prove additivity when there are no vertex connected sums, and we have that complexity is
only additive up to a certain correction summand on cyclic connected sums involving at least one
purely cyclic singular component. For this reason in this section we prove that the number
of such sums is well-defined:

\begin{prop}\label{cyc:*:prop}
Suppose that a $3$-orbifold $X$ is a connected sum of irreducible
orbifolds without vertex connected sums, and fix an efficient
realization $\rho$ of $X$ as $X_0\#\ldots\#X_n$.
For all $p\geqslant 2$ let $\nu_\rho(p)$ be the number of
$p$-cyclic connected sums in $\rho$ involving at least one singular
component without vertices.
Then $\nu_\rho(p)$ is independent of $\rho$, so it is a function of $X$ only.
\end{prop}

\begin{proof}
Two different realizations of $X$ are related by the following operations:
\begin{itemize}
\item Reordering of the $X_i$'s, without modification of the efficient system of spherical
$2$-orbifolds along which $X$ is split.
\item Modification of the splitting system according to one of the moves
$\alpha$ and $\beta$
described in Section~\ref{uniqueness:section} (see
Figg.~\ref{alphamove:fig} and~\ref{betamove:fig}).
\end{itemize}
We will show that $\nu_\rho(p)$ remains unchanged under both operations.

To deal with the first operation, we define a modified version $T_p$ of
the tree associated to the realization $X=X_0\#\ldots\#X_n$.
To construct $T_p$ we associate to each $X_i$
one node for each cyclic component of $S(X_i)$ and one
for the union of the non-cyclic components of $S(X_i)$.
We denote by $C$ the set of nodes of cyclic type
and by $N$ the set of nodes of non-cyclic type.
We now define $T_p$ by
taking an edge for each $p$-cyclic connected sum, with the
ends of the edge joining the singular components involved in the sum.

The ideas behind the construction of $T_p$ are as follows.
First, we cannot create a new $p$-cyclic
singular component by
performing an ordinary connected sum or a $q$-cyclic connected sum for
$q\neq p$.  Second, a cyclic connected sum between two cyclic components
gives one cyclic component,
while a cyclic connected sum between a non-cyclic component
and any other one gives one or two non-cyclic components.
With these facts in mind, it is very easy to see that the invariance of
$\nu_\rho(p)$ under the first move is implied by the following
graph-theoretical statement:
\emph{Let $T$ be a tree with set of nodes $C\sqcup N$.
Let $\sigma=(e_1,\ldots,e_n)$ be an ordering of the edges of $T$.
For $k=1,\ldots,n$ define $\alpha(k)$ to be $0$ if both the
ends of $e_k$ can be joined to nodes in $N$ by edges
in $e_1,\ldots,e_{k-1}$ only, and to be $1$ otherwise.
Then $\sum_{k=1}^n\alpha(k)$ is independent of $\sigma$}.

To prove this fact, we first reduce to an easier statement.
We begin by noting that the vertices in $N$ having valence more
than 1 can be blown up as in Fig.~\ref{blowup:fig}
    \begin{figure}
    \begin{center}
    \input{blowup.pstex_t}
    \nota{This move does not change the sum of the values of $\alpha$.} \label{blowup:fig}
    \end{center}
    \end{figure}
without affecting the sum in question. This allows us to assume that all
the nodes in $N$ are ``external,'' \emph{i.e.} $1$-valent. Moreover, using
induction, we can suppose that $N$ consists precisely of the external nodes,
because if an edge $e$ has an end which is external and belongs to $C$ then
$\alpha(e)=1$ whatever the position of $e$ in the ordering,
and the value of $\alpha$ on the remaining edges does not depend on
the position of $e$ in the ordering.

Switching from $\alpha(e)$ to $\beta(e)=1-\alpha(e)$ it is then sufficient
to establish the following result. The proof we present here is due
to Gw\'ena\"el Massuyeau.

\begin{lemma}\label{graph:lem}
Let $T$ be a tree. Let $E$ denote the set of
external nodes of $T$. Let $\sigma=(e_1,\ldots,e_n)$ be any
ordering of the edges of $T$. Define $\beta(k)$ to be $1$ if the ends
of $e_k$ can be joined to $E$ through edges in $e_1,\ldots,e_{k-1}$ only, and to
be $0$ otherwise. Then $\sum_{k=1}^n\beta(k)$ is independent of $\sigma$, and precisely
equal to $\#E-1$.
\end{lemma}

\begin{proof}
Attach a circle to $T$ by gluing the nodes in $E$ to $\#E$
arbitrarily chosen points on the circle, and denote by $Y$ the result.
Since $T$ is a tree, we have $\chi(Y)=1-\#E$.
Let $\sigma=(e_1,\ldots,e_n)$ be an ordering of the edges of $T$, let
$Y_k$ be the union of the circle with $e_1,\ldots,e_k$, and
let $\widetilde Y_k$ be the connected component of $Y_k$ which
contains the circle. Of course
$\widetilde Y_0$ is the circle and $\widetilde Y_k=Y$. We now have
the identity
$$ \chi(\widetilde Y_k) = \chi(\widetilde Y_{k-1}) - \beta(k) $$
which is easily established by considering separately the cases where
$e_k$ has $0$, $1$, or $2$ ends on $\widetilde Y_{k-1}$, and using
again the fact that $T$ is a tree. The identity implies that
$$\sum_{k=1}^n\beta(k)=\chi(\widetilde Y_0)-\chi(\widetilde Y_k)=\#E-1.$$
\end{proof}

To conclude the proof of Proposition~\ref{cyc:*:prop} we must now show that
$\nu_\rho(p)$ is not affected by the moves $\alpha$ and $\beta$.
Knowing the independence of the ordering,
this is actually very easy. For both the moves, we have
two realizations $\rho_0$ and $\rho_1$ of $X$ leading to
splitting systems $\calS_0$ and $\calS_1$ which
differ for one component only, so $\calS_0=\{\Sigma_0\}\cup \calS$ and
$\calS_1=\{\Sigma_1\}\cup \calS$. We then compute $\nu_{\rho_j}(p)$ by performing
first the connected sum which corresponds to $\Sigma_j$. From Figg.~\ref{alphamove:fig}
and~\ref{betamove:fig}
one sees that the connected sums along $\Sigma_0$ and $\Sigma_1$ give the same contributions
to $\nu$. After performing them, all other sums in $\rho_0$ and $\rho_1$ are the same,
which eventually proves the proposition.
\end{proof}

\begin{rem}
\emph{Proposition~\ref{cyc:*:prop} is false if one allows also
vertex connected sums in the realization of $X$, because a vertex
connected sum can create singular components without vertices.}
\end{rem}

\vspace{.5cm}

\noindent
Dipartimento di Matematica Applicata, Universit\`a di Pisa\\
Via Bonanno Pisano 25B, 56126 Pisa, Italy\\
petronio@dm.unipi.it

\end{document}